\documentclass[a4paper,12pt,reqno]{amsart}
\usepackage{amsfonts}
\usepackage{amsmath}
\usepackage{amssymb}
\usepackage[a4paper]{geometry}
\usepackage{mathrsfs}
\usepackage{csquotes}
\usepackage{enumitem}
\usepackage[colorlinks]{hyperref}
\renewcommand\eqref[1]{(\ref{#1})} %Need with hyperref
\numberwithin{equation}{section}
\theoremstyle{plain}
\newtheorem{thm}{Theorem}[section]

\theoremstyle{definition}

\newcommand{\Rn}{\mathbb R^{n}}

\begin{document}
	
	\title[Multidimensional inverse Cauchy problems for evolution equations]
	{Multidimensional inverse Cauchy problems \\ for evolution equations}
	
	\author[M. Karazym]{Mukhtar Karazym}
	\address{
		Mukhtar Karazym:
		\endgraf
		Department of Fundamental Mathematics
		\endgraf
		Gumilyov Eurasian National University, Kazakhstan
		\endgraf
		and
			\endgraf
		Department of Mathematics
		\endgraf
		Nazarbayev University, Kazakhstan
		\endgraf
		{\it E-mail address} {\rm mukhtar.karazym@nu.edu.kz}
	}
	
	\author[T. Ozawa]{Tohru Ozawa}
	\address{
		Tohru Ozawa:
		\endgraf
		Department of Applied Physics
		\endgraf
		Waseda University
		\endgraf
		Tokyo 169-8555
		\endgraf
		Japan
		\endgraf
		{\it E-mail address} {\rm txozawa@waseda.jp}
	}

	\author[D. Suragan]{Durvudkhan Suragan}
	\address{
		Durvudkhan Suragan:
		\endgraf
		Department of Mathematics
		\endgraf
		Nazarbayev University, Kazakhstan
		\endgraf
		{\it E-mail address} {\rm durvudkhan.suragan@nu.edu.kz}
	} 
	
	\thanks{The authors were supported
		by the Nazarbayev University program 091019CRP2120. No new data was collected or generated during the course of research.}
	
	
	\begin{abstract} We discuss inverse problems to finding the time-dependent coefficient for the multidimensional Cauchy problems for both strictly hyperbolic equations and polyharmonic heat equations. We also extend our techniques to the general inverse Cauchy problems for evolution equations.
	\end{abstract}
	\maketitle
	
	\tableofcontents
	
	\section{Introduction}

In inverse problems many authors have aimed to find a unique pair $(u,f)$ or $(u,a)$ from the data given by a part of solutions, where $u$ is a solution, $f$ is a source function and $a$ is a coefficient in partial differential equations. We mention only a few, see \cite{ALZ04,CR91,Isa93,KY, Mal87,Run80,Yam99} for references. Inverse problems have been considered significantly in a form of only one initial / initial-boundary problem with some additional data. In this paper, we solve multidimensional inverse problems by considering two Cauchy problems and find a unique pair, which consists of a solution and a time-dependent coefficient in an explicit form, that is, the solution is given by the Poisson integral, the time-dependent coefficient is given by the ratio of the additional data at an internal fixed point of any open bounded set in $\Rn$. 
Historically, one of the pioneering works devoted
to finding the time-dependent coefficient in the parabolic equations is a paper by Jones \cite{Jones}. In the present paper, our method is closely related to Malyshev's approach \cite{Mal87} in which the author studied inverse problems in one-dimensional degenerate parabolic equations.

The structure of this short paper is as follows. In Section \ref{sec2} we consider an inverse Cauchy problem for strictly hyperbolic equations.  The uniqueness of the solution follows from the well-posedness of the Cauchy problem, which also implies the uniqueness of the recovered coefficient. In Section \ref{sec3}, in order to solve an inverse problem for polyharmonic heat equations, first we present its fundamental solution and discuss some properties. Then by using known results from \cite{Eid64} we establish the uniqueness of the solution of Cauchy problems as a special case and solve the inverse problem in the same way as in Section \ref{sec2}. Section \ref{sec4} is devoted to a systematic study of inverse problems for general partial differential equations. Thus, we conclude that the developed method in Sections \ref{sec2} and \ref{sec3} can be applied to more general cases.

	\section{Wave equations}\label{sec2}
	In this section, we consider the (inverse) Cauchy problem for the strictly hyperbolic equation in the form
	\begin{equation}\label{Cauchy3.1}
	\left\{ \begin{split} 
	&L_{\Phi}u:=\partial^{2}_{t}u(t,x)-\Phi(t)\Delta_{x}u(t,x)=0,\quad 0<t\leq T< \infty, \quad x\in\Rn,\\
	&u(t,x)|_{t=0}=u_{0}(x),\quad x\in\Rn,\\
	&\partial_{t} u(t,x)|_{t=0}=u_{1}(x),\quad x\in\Rn,
	\end{split}
	\right.
	\end{equation} 
	where $\Phi(t)\geq C>0$ and is Lipschitz continuous. Here and after $\Delta_{x}=\sum_{j=1}^{n}\partial^{2}_{x_{j}}$ is the Laplacian. Problem \eqref{Cauchy3.1} is well-posed in Sobolev spaces. It is known that for any $u_{0}\in H^{m}(\Rn)$ and $u_{1}\in H^{m-1}(\Rn)$ there exists a unique solution
	$u\in C([0,\,T]; H^{m}(\Rn))\cap C^{1}([0,\,T]; H^{m-1}(\Rn))$ for $m\geq 1$, see \cite[Ch. 9]{Hor63}.
	The solution of problem \eqref{Cauchy3.1} is given by
	\begin{equation}\label{solution3.1}
	u(t,x)=\int_{\Rn}E(t,0,x,y)u_{1}(y)\,dy+\int_{\Rn}\partial_{\tau}E(t,\tau,x,y)|_{\tau=0}^{\tau=T}u_{0}(y)\,dy,
	\end{equation}
	where $E$ is the fundamental solution. 
	
	Let $\Omega\subset\Rn,\,n\in\mathbb{N},$ be an open bounded set with piecewise smooth boundary $\partial\Omega$ and let $q\in\Omega$ be a fixed point throughout the paper.  Now we consider the inverse problem by studying two Cauchy problems with the additional data at $q\in\Omega$:
	\begin{equation}\label{Inverse3.1}
	\left\{ \begin{split} 
	&\partial^{2}_{t}\,v(t,x)-\Phi(t)\Delta_{x}v(t,x)=0,\quad 0<t\leq T, \quad x\in\Rn,\\
	&v(t,x)|_{t=0}=0,\quad  x\in\Rn,\\
	&\partial_{t} v(t,x)|_{t=0}=u_{1}(x),\quad x\in\Rn,\\
	&v(t,x)|_{x=q}=h_{1}(t),\quad 0< t \leq T,\;q\in\Omega,
	\end{split}
	\right.
	\end{equation}
	and
	\begin{equation}\label{Inverse3.2}
	\left\{ \begin{split} 
	&\partial^{2}_{t}\,w(t,x)-\Phi(t)\Delta_{x}w(t,x)=0,\quad 0<t\leq T, \quad x\in\Rn,\\
	&w(t,x)|_{t=0}=0,\quad  x\in\Rn,\\
	&\partial_{t} w(t,x)|_{t=0}=\Delta_{x}u_{1}(x),\quad x\in\Rn,\\
	&w(t,x)|_{x=q}=h_{2}(t),\quad 0< t \leq T,\;q\in\Omega,
	\end{split}
	\right.
	\end{equation}	
	where $u_{1}\in C^{2}(\Omega)$ and $\operatorname{supp}(u_{1})\subset \Omega$. Assuming that $\Phi$ is a certain Lipschitz continuous function and satisfies $\Phi(t)\geq C>0$, we solve Cauchy problems in \eqref{Inverse3.1} and \eqref{Inverse3.2}. Their solutions are expressed by
	\begin{equation*}
	v(t,x)=\int_{\Rn}E(t,0,x,y)(\chi_{\Omega}u_{1})(y)dy=\int_{\Omega}E(t,0,x,y)\,u_{1}(y)dy,
	\end{equation*}
	and
	\begin{equation*}
	w(t,x)=\int_{\Rn}E(t,0,x,y)\,(\chi_{\Omega}\Delta_{y} u_{1})(y)dy=\int_{\Omega}E(t,0,x,y)\,\Delta_{y} u_{1}(y)dy,
	\end{equation*}
correspondingly, where $\chi_{\Omega}$ is the characteristic function of $\Omega$. 
From the additional data it follows that
\begin{equation}\label{h1}
	v(t,q)=\int_{\Omega}E(t,0,x,y)|_{x=q}u_{1}(y)dy=h_{1}(t),
\end{equation}

	\begin{equation*}
	w(t,q)=\int_{\Omega}E(t,0,x,y)|_{x=q}\Delta_{y} u_{1}(y)dy=h_{2}(t).
	\end{equation*}
	Assuming that $h_{1}\in C^{2}[0,\,T]$ and  twice differentiating \eqref{h1}, we arrive at 
	\begin{equation*}
	h^{\prime\prime}_{1}(t)=\int_{\Omega}	\partial^{2}_{t}E(t,0,x,y)|_{x=q}\,u_{1}(y)\,dy.
	\end{equation*}
Now by using the fact that $E$ is the fundamental solution, that is, 
	\begin{equation*}
	\partial^{2}_{t}\,E(t,0,x,y)=\Phi(t)\Delta_{x}E(t,0,x,y)=\Phi(t)\Delta_{y}E(t,0,x,y),\quad t>0,
	\end{equation*}
	and also applying Green's second identity, we arrive at
	\begin{multline*}
	h_{1}^{\prime\prime}(t)=\int_{\Omega}	\partial^{2}_{t}E(t,0,x,y)|_{x=q}u_{1}(y)dy=\Phi(t)\int_{\Omega}	\Delta_{y} E(t,0,x,y)|_{x=q}u_{1}(y)dy\\=\Phi(t)\int_{\Omega}	 E(t,0,x,y)|_{x=q}\Delta_{y}u_{1}(y)dy=\Phi(t)h_{2}(t).
	\end{multline*}
	As we have assumed that $\Phi$ is Lipschitz continuous, we require
	$h_{2}$ to be a Lipschitz continuous function such that $h_{2}(t)\neq 0$ for all $t\in (0,\,T]$. As we have assumed that $\Phi(t)\geq C >0$, we require $h_{1}$ and $h_{2}$ to satisfy $\frac{h_{1}^{\prime\prime}(t)}{h_{2}(t)}\geq C>0$ for all $t\in (0,\,T]$. So, we have recovered the coefficient $\Phi$. The uniqueness of the solution of the inverse problem follows from the uniqueness of the solution of the Cauchy problems. Thus, we obtain the following theorem.
	\begin{thm}\label{thm3.1}
		Let us make the assumptions:
		\begin{enumerate}
			\item[$1)$]	$u_{1}\in C^{2}(\Omega)$ and $\operatorname{supp}(u_{1})\subset \Omega$;
			\item[$2)$] $h_{1}\in C^{2}[0,\,T]$;
			\item[$3)$]	$h_{2}$ is a  Lipschitz continuous function such that $h_{2}(t)\neq 0$ for all $t\in (0,\,T]$;
			\item[$4)$] $\frac{h_{1}^{\prime\prime}(t)}{h_{2}(t)}\geq C>0$ for all $t\in (0,\,T]$.
		\end{enumerate}	
		Then there exists a unique solution of inverse problem \eqref{Inverse3.1}-\eqref{Inverse3.2} with the corresponding Lipschitz continuous coefficient $\Phi(t)=\frac{h_{1}^{\prime\prime}(t)}{h_{2}(t)}$ for all $t\in (0,\,T]$.
	\end{thm}

	\section{Polyharmonic heat equations}\label{sec3}
	The idea of the proof from the previous section can be applied to inverse Cauchy problems for more general evolution equations. To demonstrate it in this section, we discuss an inverse problem for the Cauchy problem for polyharmonic heat equations with a time-dependent coefficient.
	
	 We consider the following Cauchy problem for the polyharmonic heat equation with the time-dependent coefficient
	\begin{equation}\label{Cauchy1}
	\left\{ \begin{split} 
	&\partial_{t}u(t,x)+\alpha(t)\big(-\Delta_{x}\big)^{m}u(t,x)=0,\quad 0<t\leq T, \quad x\in\Rn,\\
	&u(t,x)|_{t=0}=u_{0}(x),\quad x\in\Rn,
	\end{split}
	\right.
	\end{equation} 
	where $m\in \mathbb{N}$, $u_{0}$ is a given function. Here the coefficient $\alpha\in L^{1}[0,T]$ satisfies the assumption that
	\begin{equation*}
	\alpha_1(t):=\int_{0}^{t}\alpha(s)ds>0,\quad \text{for all}~0<t \leq T.
	\end{equation*}
	By the Fourier transfrom, see \cite{RE66}, one finds its fundamental solution for all $t>0$ and $x\in\Rn$ in the form 
	\begin{equation}\label{fundamental solution}
	E_{\alpha_{1}}(t,x):=(2\pi)^{-n}\int_{\Rn}e^{ix\cdot s-|s|^{2m}\alpha_{1}(t)}ds.
	\end{equation}
	Note that the fundamental solution \eqref{fundamental solution} can be reduced to the one-dimensional integral representation
	\begin{equation}\label{fundamental solution1}
	E_{\alpha_{1}}(t,x)=(2\pi)^{-\frac{n}{2}}\alpha_{1}(t)^{-\frac{n}{2m}}\int_{0}^{\infty}e^{-r^{2m}}r^{\frac{n}{2}}\,J_{\frac{n-2}{2}}\big(r\,|x|\,\alpha_{1}(t)^{-\frac{1}{2m}}\big)\,dr,
	\end{equation}
	where $J_{k}$ is the Bessel function of the first kind. For details, we refer to \cite[p. 183-184]{EZ98}. From \eqref{fundamental solution1} it is obvious that $E_{\alpha_{1}}(t,x-y)=E_{\alpha_{1}}(t,y-x)$ for all $x,\,y\in \Rn$ and $t>0$.
	
	In \cite{Mal91}, the author studied Cauchy problem \eqref{Cauchy1} in the case when $m=n=1$. In our previous paper \cite{KS19}, we investigated Cauchy problem \eqref{Cauchy1} in the case when $m=1$, $n\geq 2$.
	
	Now we assume that the equation in \eqref{Cauchy1} is uniformly parabolic in the sense of Petrovskii, see \cite{Eid64,EZ98} for a precise definition. We recall the following theorem, which is a particular case of \cite[Theorem 5.3]{Eid64}.
	\begin{thm}\cite{Eid64}\label{thm2.1}
	Let $\alpha\in C[0,\,T]$ and $u_{0}\in C^{2m,\,\gamma}(\Omega),\,0<\gamma<1$, with compact support in $\Omega$. Then there exists a unique solution $u\in C^{2m,\,\gamma,\,0}([0,\,T]\times\Omega)$ of problem \eqref{Cauchy1} defined by the Poisson integral 
	\begin{equation}\label{solution1}
	u(t,x)=\int_{\Rn} E_{\alpha_{1}}(t,x-y)(\chi_{\Omega}u_{0})(y)dy=\int_{\Omega} E_{\alpha_{1}}(t,x-y)u_{0}(y)dy.
	\end{equation}
	\end{thm}
	Now we consider an inverse problem for \eqref{Cauchy1}. Let us consider Cauchy problem \eqref{Cauchy1} with the additional data 
	\begin{equation}\label{data1}
	u(t,x)|_{x=q}=h_{1}(t),\quad 0\leq t \leq T,\; q\in \Omega, 
	\end{equation}
	and another Cauchy problem
	\begin{equation}\label{Cauchy2}
	\left\{ \begin{split} 
	&\partial_{t}v(t,x)+\alpha(t)\big(-\Delta_{x}\big)^{m}v(t,x)=0,\quad 0<t\leq T, \quad x\in\Rn,\\
	&v(t,x)|_{t=0}=-\big(-\Delta_{x}\big)^{m}u_{0}(x),\quad x\in\Rn,
	\end{split}
	\right.
	\end{equation} 
	with the additional data 
	\begin{equation}\label{data2}
	v(t,x)|_{x=q}=h_{2}(t),\quad 0\leq t \leq T,\;q\in\Omega.
	\end{equation}
	
	Let us suppose that $\alpha\in C[0,\,T]$ which ensures uniform parabolicity of the homogeneous equation in the sense of Petrovskii. In addition, we assume that $u_{0}\in C^{4m,\,\gamma}(\Omega)$ with $\operatorname{supp}(u_{0})\subset \Omega$. Then, we have unique solutions of Cauchy problems \eqref{Cauchy1}, \eqref{Cauchy2} and they can be represented by formula \eqref{solution1} and by
	\begin{equation*}
	v(t,x)=-\int_{\Rn} E_{\alpha_{1}}(t,\,x-y)\,(\chi_{\Omega}	\left(-\Delta_{y}\right)^{m}u_{0})(y)\,dy=-\int_{\Omega} E_{\alpha_{1}}(t,\,x-y)\,\left(-\Delta_{y}\right)^{m}u_{0}(y)\,dy.
	\end{equation*}
	Using the additional data \eqref{data1}, \eqref{data2}, we come to
	$$u(t,x)|_{x=q}=\int_{\Omega} E_{\alpha_{1}}(t,x-y)|_{x=q}\,u_{0}(y)dy=h_{1}(t),$$
	and
	$$v(t,x)|_{x=q}=-\int_{\Omega} E_{\alpha_{1}}(t,x-y)|_{x=q}\,\big(-\Delta_{y}\big)^{m}u_{0}(y)dy=h_{2}(t).$$
	Now we assume that $h_{1}\in C^{1}[0,\,T]$.
	Differentiating $h_{1}$, and then using the fact that the fundamental solution is an even function with respect to spatial variables and Green's second identity, we arrive at
	\begin{eqnarray}\label{relation}
	\begin{aligned}
	&h_{1}^{\prime}(t)=\int_{\Omega}\partial_{t} E_\alpha(t,x-y)|_{x=q} \,u_{0}(y) dy=-\alpha(t)\int_{\Omega}\big(-\Delta_{y}\big)^{m} E_\alpha(t,x-y)|_{x=q}u_{0}(y)  dy
	\\
	&=-\alpha(t)\int_{\Omega}\big(-\Delta_{y}\big)^{m-1} E_\alpha(t,x-y)|_{x=q}(-\Delta_{y}\big)u_{0}(y)dy
	\\
	&=...=-\alpha(t)\int_{\Omega} E_\alpha(t,x-y)|_{x=q} \big(-\Delta_{y}\big)^{m}u_{0}(y)dy=\alpha(t)h_{2}(t).
	\end{aligned}
	\end{eqnarray}
	Here we assume that $h_{2}\in C[0,\,T]$ and $h_{2}(t)\neq 0$ for all $0\leq t \leq T$. Thus, we have
	
	\begin{equation}\label{coef}
	\alpha(t)=\frac{h_{1}^{\prime}(t)}{h_{2}(t)},\quad 0\leq t\leq T.
	\end{equation}
	
	Note that $h_{2}(0)=v(0,x)|_{x=q}=-\big(-\Delta_{x}\big)^{m}u_{0}(x)|_{x=q}\neq 0$.
	By the assumption, $\alpha$ ensures that the equation in \eqref{Cauchy1} is uniformly parabolic in the sense of Petrovskii. The uniqueness of the inverse problem follows from the uniqueness of $\alpha$. That is, we have proved the following theorem.
	
	\begin{thm}\label{thm2.2}
	Let us make the following assumptions:
	\begin{enumerate}
		\item[$1)$] $u_{0}\in C^{4m,\,\gamma}(\Omega)$ and $\operatorname{supp}(u_{0})\subset\Omega$;
		\item[$2)$] $h_{1}\in C^{1}[0,\,T]$;
		\item[$3)$] $h_{2}\in C[0,\,T]$ such that $h_{2}(t)\neq 0$ for all $0\leq t \leq T$ (which also implies $v(0,x)|_{x=q}=-\big(-\Delta_{x}\big)^{m}u_{0}(x)|_{x=q}=h_{2}(0)\neq 0$);
		\item[$4)$] $\frac{h_{1}^{\prime}}{h_{2}}$ ensures that the equation in \eqref{Cauchy1} is uniformly parabolic in the sense of Petrovskii.
	\end{enumerate}		
	Then there exists a unique solution for inverse problem \eqref{Cauchy1}, \eqref{data1}-\eqref{data2} with the coefficient $\alpha\in C[0,\,T]$ defined by \eqref{coef}.
	\end{thm}	
\section{General Cauchy problems}\label{sec4}
Let us consider the higher order linear partial differential equation with the first order partial time derivative
\begin{equation}\label{hp equation}
L_{\Psi}u(t,x):=\partial_{t}u(t,x)-\Psi(t)L_{x}[u](t,x)=0,\quad t>0,\quad x\in \Rn,
\end{equation}
where 
$$L_{x}[u](t,x):=\sum_{|k|\leq m} A_{k}(t,x)\,\partial^{k} u(t,x),$$
in multi-index notation with
$k=(k_{1},\ldots,k_{n})\in (\mathbb{N}\cup \{0\})^{n}, \;  |k|=k_{1}+\ldots+k_{n},\; \partial^{k}= \partial_{1}^{k_{1}}\ldots \partial_{n}^{k_{n}}$ and $\partial_{j}=\frac{\partial}{\partial x_{j}}$.
We assume that the coefficients $\Psi(t)A_{k}(t,x)$ are sufficiently smooth and bounded if necessary. Also we assume that the coefficients of the highest derivatives are nonzero everywhere in $\Rn$.
The solution of \eqref{hp equation} with the initial condition
\begin{equation}\label{hp initial}
u(t,x)|_{t=0}=u_{0}(x),\quad x\in\Rn,
\end{equation}
is given by
\begin{equation*}
u(t,x)=\int_{\Rn}\varepsilon_{1}(t,0,x,y)u_{0}(y)dy,
\end{equation*}
where $\varepsilon_{1}$ is the fundamental solution, see \cite[Section 9.6.3-1]{Pol02}. We assume that the solution of \eqref{hp equation}-\eqref{hp initial} is unique and the fundamental solution 
 satisfies
\begin{equation*}
\varepsilon_{1}(t,0,x,y)=\varepsilon_{1}(t,0,y,x),\quad t>0.
\end{equation*}
Note that, the fundamental solution is a solution of
\begin{equation*}
\left\{ \begin{split} 
&\partial_{t}\varepsilon_{1}(t,\tau,x,y)-\Psi(t)L_{x}[\varepsilon_{1}](t,\tau,x,y)=0,\quad t>\tau\geq 0, \quad x,\,y\in\Rn,\\
&\varepsilon_{1}(t,\tau,x,y)|_{t=\tau}=\delta(x-y),\quad x,\,y\in\Rn.
\end{split}
\right.
\end{equation*}
We also assume that the operator $L_{x}$ satisfies Green's second identity, that is,
\begin{equation*}
\int_{\Omega} (\varphi L_{x}[\varepsilon_{1}]-\varepsilon_{1} L_{x}[\varphi])dx=0,
\end{equation*}
for all $\varphi\in C^{|k|}(\Omega)$ and $\operatorname{supp}(\varphi)\subset \Omega$. We solve the following inverse problem of finding a unique pair $(u,\Psi)$
\begin{equation}\label{Inverse4.1}
\left\{ \begin{split} 
&\partial_{t}u(t,x)-\Psi(t)L_{x}[u](t,x)=0,\quad 0<t\leq T, \quad x\in\Rn,\\
&u(t,x)|_{t=0}=u_{0}(x),\quad x\in\Rn,\\
&u(t,x)|_{x=q}=h_{1}(t),\quad 0\leq t \leq T,\;q\in \Omega,
\end{split}
\right.
\end{equation}
and
\begin{equation}\label{Inverse4.2}
\left\{ \begin{split} 
&\partial_{t}v(t,x)-\Psi(t)L_{x}[v](t,x)=0,\quad 0<t\leq T, \quad x\in\Rn,\\
&v(t,x)|_{t=0}=L_{x}[u_{0}](x),\quad x\in\Rn,\\
&v(t,x)|_{x=q}=h_{2}(t),\quad 0\leq t \leq T,\;q\in \Omega,
\end{split}
\right.
\end{equation}
where $u_{0}\in C^{|k|}(\Omega)$ with $\operatorname{supp}(u_{0})\subset \Omega$ and $h_{1},\,h_{2}$ are sufficiently smooth such that $h_{2}(t)\neq 0$ for all $t\in [0,\,T]$ which implies $h_{2}(0)=L_{x}[u_{0}](x)|_{x=q}\neq 0$.
\begin{thm}\label{thm4.1}
	Let us assume
	\begin{enumerate}
		\item[$1)$] $L_{\Psi}$, $\varepsilon_{1}$ and $L_{x}$ satisfy all the previous assumptions;
		\item[$2)$]  $u_{0}\in C^{|k|}(\Omega)$ with $\operatorname{supp} \,u_{0}\subset \Omega$;
		\item[$3)$] $h_{1},\,h_{2}$ are sufficiently smooth such that $h_{2}(t)\neq 0$ for all $t\in [0,\,T]$.
	\end{enumerate}
Then there exists a unique solution of inverse problem \eqref{Inverse4.1}-\eqref{Inverse4.2} with the corresponding coefficient $\Psi(t)=\frac{h^{\prime}(t)}{h_{2}(t)}$.
	\end{thm}

Also one can consider an inverse problem for higher order equations with the second order partial time derivative, see \cite[Section 9.6.3-2]{Pol02}, that is,
\begin{equation}\label{Inverse4.3}
\left\{ \begin{split} 
&L_{\Lambda}u(t,x):=\partial^{2}_{t}\,u(t,x)-\Lambda(t)L_{x}[u](t,x)=0,\quad 0<t\leq T, \quad x\in\Rn,\\
&u(t,x)|_{t=0}=0,\quad x\in\Rn,\\
&\partial_{t}u(t,x)|_{t=0}=u_{1}(x),\quad x\in\Rn,\\
&u(t,x)|_{x=q}=h_{1}(t),\quad 0< t \leq T, \;q\in \Omega,
\end{split}
\right.
\end{equation}
and 
\begin{equation}\label{Inverse4.4}
\left\{ \begin{split} 
&\partial^{2}_{t}v(t,x)-\Lambda(t)L_{x}[v](t,x)=0,\quad 0<t\leq T, \quad x\in\Rn,\\
&v(t,x)|_{t=0}=0,\quad x\in\Rn,\\
&\partial_{t}v(t,x)|_{t=0}=L_{x}[u_{1}](x),\quad x\in\Rn,
\\
&v(t,x)|_{x=q}=h_{2}(t),\quad 0< t \leq T,\;q\in \Omega,
\end{split}
\right.
\end{equation}
where the coefficients $\Lambda(t)A_{k}(t,x)$ are sufficently smooth and bounded if necessary, $u_{1}\in C^{|k|}(\Omega)$ with $\operatorname{supp}(u_{1})\subset \Omega$.
Note that, the fundamental solution $\varepsilon_{2}(t,\tau,x,y)$ of \eqref{Inverse4.3}-\eqref{Inverse4.4} solves 
\begin{equation*}
\left\{ \begin{split} 
&\partial^{2}_{t}\varepsilon_{2}(t,\tau,x,y)-\Lambda(t)L_{x}[\varepsilon_{2}](t,\tau,x,y)=0,\quad t>\tau\geq 0, \quad x,\,y\in\Rn,\\
&\varepsilon_{2}(t,\tau,x,y)|_{t=\tau}=0,\quad x,\,y\in\Rn,\\
&\partial_{t}\varepsilon_{2}(t,\tau,x,y)|_{t=\tau}=\delta(x-y),\quad x,\,y\in\Rn.
\end{split}
\right.
\end{equation*}
Assuming all necessary facts for $L_{\Lambda}$, $L_{x}$ and $\varepsilon_{2}$, we can obtain a unique solution of \eqref{Inverse4.3}-\eqref{Inverse4.4} with the corresponding time-dependent sufficiently smooth coefficient $\Lambda(t)=\frac{h_{1}^{\prime\prime}(t)}{h_{2}(t)}.$

Thus the above idea can be extended to general inverse Cauchy problems.

\section{Conclusion}

We propose a method to finding the time-dependent coefficient for multidimensional Cauchy problems for both strictly hyperbolic equations and polyharmonic (degenerate) heat equations. In addition, we show that the method can be extended to solving the general multidimensional inverse Cauchy problems for higher order linear partial differential (evolution) equations.   Basically, it is a theoretical work. However, it can be also applied to some concrete real-world processes since we consider a general class of the inverse Cauchy problems for evolution equations. The key idea is to apply the Poisson integral to find a unique pair that consists of a solution and a time-dependent coefficient in an explicit form. The time-dependent coefficient is given by the ratio of the additional data at an internal fixed point of any open bounded set of the Euclidean space. We believe a new numerical algorithm to finding the time-dependent coefficient for the Cauchy problems for evolution equations can be tested by the explicit solutions given in this paper.


\begin{thebibliography}{HOHOLT08}
	\bibitem{ALZ04}
	H. Azari, Ch. Li, Y. Nie, S. Zhang, Determination of an unknown coefficient in a
	parabolic inverse problem, {\em Dyn. Contin. Discrete Impuls. Syst., Ser. A, Math. Anal.} 11
	(2004), 665--674.
	
	\bibitem{CR91}
	J. R. Cannon, W. Rundell, Recovering a time dependent coefficient in a parabolic differential
	equation, {\em J. Math. Anal. Appl.} 160 (1991), 572--582.
	
	\bibitem{Eid64}
	S. D. Eidelman, Parabolic Systems, Nauka, Moscow, 1964 (in Russian).
			
	\bibitem{EZ98}
	S. D. Eidelman, N. V. Zhitarashu, Parabolic Boundary Value Problems, Oper. Theory
	Adv. Appl., 101, Birkhauser, Basel, 1998.
		
	\bibitem{Hor63}
	L. H\"ormander, Linear Partial Differential Operators,  Springer-Verlag, Berlin, 1963.
	
	\bibitem{Isa93}
	V. Isakov,  Uniqueness and stability in multi-dimensional inverse problems, {\em Inverse Problems}, 9 (1993): 579--621. 
		
		
		\bibitem{Jones}	B. F. Jones, The determination of a coefficient in a parabolic differential equation. Part I: Existence and uniqueness,
		{\em J. Math. Mech}, 11 (1962), 907--918.
		
		
	\bibitem{KS19}
	M. Karazym, D. Suragan, Trace formulae of potentials for degenerate parabolic equations, pre-print, 2020.

	\bibitem{KY}
	Y. Kian, M. Yamamoto. 
	Reconstruction and stable recovery of source terms and coefficients appearing in diffusion equations.
	{\em Inverse Problems}, 35(11) (2019):1--14.

	\bibitem{Mal87}
	I. Malyshev.
	On the inverse problem for a heat-like equation, {\em J. of Appl. Math. and Simulation},  1(2) (1987):81--97.	

	\bibitem{Mal91}
	I. Malyshev,
	On the parabolic potentials in degenerate-type heat equations.
	{\em J. Appl. Math. Stoch. Anal.}, 2(4) (1991):147--160.	
	
	\bibitem{Pol02}
	A. D. Polyanin, Handbook of Linear Partial Differential Equations for Engineers and Scientists.
	Chapman  Hall/CRC Press, 2002.
	
	\bibitem{RE66}
	V. D. Repnikov, S. D. Eidel’man, Necessary and 
	sufficient conditions for establishing a solution to
	the Cauchy problem, {\em Dokl. Akad. Nauk SSSR}, 1966,
	Volume 167, Number 2, 298--301 (in Russian).


	\bibitem{Run80}
	W. Rundell, Determination of an unknown nonhomogeneous term in a linear partial differential equation from overspecified boundary data, {\em Appl. Anal.}, 10 (1980), 231--242.
	
	\bibitem{Yam99}
	M. Yamamoto, Uniqueness and stability in multidimensional hyperbolic inverse
	problems, {\em J. Math. Pures Appl.} 78 (1999), 65--98.
	
	
	\end{thebibliography}
\end{document}